\documentclass{article}
\usepackage{amsmath,amssymb,amsfonts}
\newtheorem{theorem}{Theorem}
\newtheorem{lemma}{Lemma}[theorem]
\newtheorem{remark}{Remark}
\usepackage{algpseudocode}
\usepackage{algorithm}
\usepackage{geometry,graphicx}
\usepackage{textcomp}
\usepackage{nicefrac}
\usepackage{xcolor}
\usepackage{cancel}
\usepackage{url}
\usepackage{bm}
\floatname{algorithm}{Algorithm}
\title{Optimal Control of H-Mode Tokamak Plasma Temperature based on Pontryagin's Principle}
\author{Slim Jmal, Matteo Tacchi-Bénard and Emmanuel Witrant}

\begin{document}

\maketitle

\begin{abstract}
This paper studies the decay of an objective functional using a new control technique within Pontryagin's framework. Convergence analysis is carried out on the infinite-dimensional space of Tokamak plasma dynamical state as described by weakly decoupled nonlinear partial differential equations. An adjoint-based optimal control is derived to minimize the deviation from a predefined dynamical trajectory leading to the desired target state at stationary regime, by turning Pontryagin's transversality conditions into a continuum of horizons. A feedback controller is proposed to steer the system efficiently in real time, as opposed to an open-loop controller resulting from the classical Pontryagin's setting. An algorithm synthesizing the constraint-free optimal controller is used for profile tracking based on experimental data.

\begin{center}
    \textbf{Keywords}
\end{center}
Thermonuclear fusion, nonlinear PDEs, H-mode Tokamak,  Pontryagin's principle.

\begin{center}
    \textbf{Acknowledgements}
\end{center}
This work has been carried out within the framework of the EUROfusion Consortium, funded by the European Union via the Euratom Research and Training Programme (Grant Agreement No 101052200 - EUROfusion). Views and opinions expressed are however those of the author(s) only and do not necessarily reflect those of the European Union or the European Commission. Neither the European Union nor the European Commission can be held responsible for them. The work of Slim Jmal is supported by the program "Initiatives de Recherche à Grenoble Alpes" through the grant SOS-Fusion.
\end{abstract}

\section{Introduction}
Controlling electronic temperature of fusion plasma inside Tokamaks is crucial for enhancing their performance and achieving breakeven for future commercial viability. High-confinement modes (H-mode) are particularly characterized by higher energy at the plasma center, leading to significantly improved confinement times as \cite{keilhacker1987h} showed. Maintaining a stable temperature profile in this regime is critical to prevent instabilities while ensuring efficient heating and transport control.

Different strategies have been explored for controlling plasma profiles. Early works by \cite{moreau2008two} focused on linearized models using singular perturbation theory to separate time scales and control the magnetic and kinetic profiles. Other approaches like that of \cite{moreau2013integrated, witrant2007control} relied on first-principles models that capture the dominant nonlinear dynamics of the plasma. To account for the distributed nature of plasma transport, \cite{felici2011real} among others used spatial discretization methods for profile control, whereas \cite{boyer2013first,argomedo2013lyapunov} employed infinite-dimensional PDE-based control techniques.

While previous studies, for example those of \cite{mavkov2017distributed, mavkov2018experimental}, have addressed plasma dynamics coupling through model reduction and data-driven approximations, our focus is on nonlinear PDE-based control, specifically within \cite{pontryagin1962mathematical}'s optimal control framework. This provides a systematic way to derive optimal control laws by solving a state-adjoint system of constrained PDEs. Instead of using a single target profile, we introduce the concept of a continuum of horizons, dynamically adapting the reference trajectory over time to improve stability and tracking performance.

Control of parabolic partial differential equations (PDEs) in finite time is a central problem in mathematical control theory, with wide-ranging applications spanning from physics to engineering. Among the major methods for achieving that, Pontryagin Maximum Principle (PMP) has emerged as a fundamental tool, extending beyond the finite-dimensional control of ordinary differential equations (ODEs). Indeed, foundational developments through direct and variational methods were made by \cite{fursikov1996controllability}, where they relied on Hilbert spaces for exact and approximate controllability. Further advances concerning optimal control of distributed parameter systems are presented by \cite{li1991}, offering necessary conditions for control in infinite-dimensional settings.

In the specific case of diffusion-type equations, applying PMP involves the coupling of the state dynamics with its adjoint state equation evolving backwards, corresponding to Lagrange multipliers associated with PDE constraints, as rigorously analyzed by \cite{barbu1993}. For semilinear parabolic systems with distributed or boundary controls, a detailed treatment of necessary optimality conditions and numerical methods was done by \cite{troltzsch2010optimal}. In addition, the seminal work of \cite{lions1971optimal} established the abstract infinite-dimensional PMP over Hilbert spaces, providing a foundational framework that continues to influence modern developments, including extensions to sparsity and pointwise control constraints, as discussed by \cite{casas2012sparse}.

More recent developments have been made to enrich the theory of \cite{pontryagin1962mathematical}. Time optimal control problems for abstract parabolic systems with perturbations were analyzed by \cite{tucsnak2016perturbations} leveraging PMP, highlighting the sensitivity of optimal trajectories and control structures under parameter variations while providing key insights on the robustness of his formulation. Furthermore, the work of \cite{bonnans2021state} provides a refined PMP-based analysis with second-order optimality conditions for the state-constrained control of semilinear parabolic equations.

Model Predictive Control (MPC) strategies for parabolic PDEs have recently gained parallel attention for their ability to iteratively solve optimal control problems with a receding target, as studied by \cite{dubljevic2006predictive}. Receding Horizon Control (RHC) schemes for nonlinear parabolic PDEs with boundary control inputs, specifically, were examined by \cite{hashimoto2012receding}. More recently, the stabilizing properties of receding horizon strategies in infinite-dimensional settings, particularly for reaction-diffusion systems, have been analyzed by \cite{ito2002receding}, where control Lyapunov functionals and terminal costs are incorporated to ensure convergence and stability.

Building on these latter advancements in the foundational theory of \cite{pontryagin1962mathematical}, this work proposes a new technique for solving optimal control problems of nonlinear parabolic PDEs, particularly the diffusion equation of Tokamaks' plasma temperature, integrating MPC-like feedback mechanisms by continuously tracking reference target states. While retaining the structure of PMP in infinite-dimensional Hilbert spaces, this strategy blends the more recent techniques from RHC. As it is defined, our approach is particularly suited for scenarios with functional-norm control constraints, and even geometric shape constraints but with further technical adaptations.

The paper is structured as follows: Section \ref{sec:setting} describes the plasma temperature dynamics and formulates the optimal control problem. Section \ref{sec:formalism} introduces our Pontryagin-based approach for solving the problem as well as a regularity/boundedness analysis and a convergence analysis. Section \ref{sec:results} presents the algorithm and simulation results. Section \ref{sec:concl} concludes with future research directions.
\section{System Description and Control Problem} \label{sec:setting}
\subsection{Electron Temperature Dynamics}
The evolution of the electron temperature $T_e(x,t)$ is governed by a nonlinear parabolic PDE modeling heat transport in Tokamak plasmas. Under the assumption of toroidal symmetry, the system reduces to the 1D radial diffusion equation studied by \cite{clemenccon2004analytical} in the cylindrical coordinate system with radial variable $x \in [0,1]$
\begin{equation} \label{eq:dynamics}
    \frac{3}{2} \frac{\partial (n_e T_e)}{\partial t} = \frac{1}{a^2} \frac{1}{x} \frac{\partial}{\partial x} \biggl( x n_e \chi_e \frac{\partial T_e}{\partial x} \biggr) - P_{sink} + P_{sources}
\end{equation}
with mixed Dirichlet and Neumann boundary conditions $T_e(1,t) = T_{edge}(t) \approx 0, \frac{\partial T_e}{\partial x}(0,t) = 0, \, \forall t \geq 0$ and
\begin{enumerate}
    \item[\textbullet] $a$ is the minor radius of the vacuum chamber.
    \item[\textbullet] $n_e(x,t)$ is the electron density profile.
    \item[\textbullet] $\chi_e(x,t)$ represents the electron heat diffusivity.
    \item[\textbullet] $P_{sources}(x,t)$ includes external heating sources such as Ohmic heating ($P_{OH}$), auxiliary power ($P_{aux}$), and Neutral Beam Injection.
    \item[\textbullet] $P_{sink}(x,t)$ accounts for energy loss mechanisms, such as electron-ion equipartition losses and radiative cooling, which are often neglected in simplified models used by \cite{witrant2007control} \cite{felici2011real}.
\end{enumerate}
\subsection{Electron Heat Diffusivity Model} \label{subsec:model}
Semi-empirical models are employed instead of a fully analytic model for the electron heat diffusivity because of the complexity of plasma heat transport. In this work, we adopt an extended Bohm/gyro‐Bohm model developed by \cite{christofides2002nonlinear, pianroj2012simulations} that has been successfully used in transport simulations of H-mode Tokamak plasmas. The model expresses the electron heat diffusivity as
\begin{equation} 
    \chi_e = \chi_{ec} \times f_s, \quad \chi_{ec} = ( 2 \chi_{Be} + \chi_{gBe} ) f_s.\end{equation} 
where the classical diffusivity $\chi_{e_c}$ is decomposed into Bohm and gyro‐Bohm contributions. The latter writes
\begin{equation} 
    \chi_{gBe} = 5 \times 10^{-6} \sqrt{T_e} \biggl| \frac{\nabla T_e}{B_{\phi_0}^2} \biggr|
\end{equation}
whereas the Bohm diffusivity is expressed as
\begin{equation} 
    \chi_{Be} = 4 \times 10^{-5} R \biggl| \frac{\nabla (n_e T_e)}{n_e B_{\phi_0}} \biggr| q^2 \biggl(\frac{T_{e,0.8}-T_{e,1}}{T_{e,1}} \biggr) 
\end{equation}
where $B_{\phi_0}$ is the toroidal magnetic field, $R$ is the major radius, $q$ is the safety factor and $T_{e,1}$ (resp $T_{e,0.8})$ represents the electron temperature at $x=1$ (resp $0.8)$ and the last ratio represents the phenomena in which the diffusivity decreases when the edge temperature is increased. \\ \\
The suppression function $f_s(x)$ accounts for the reduced transport due to turbulence stabilization mechanisms
\begin{equation} 
    f_s(x) = \frac{1}{1 + k \left( \frac{\omega_{E\times B}}{\gamma_{ITG}} \right)^2} \times \frac{1}{\max \left( 1,(s-s_{thres})^2 \right)}
\end{equation} 
where
\begin{itemize} 
    \item[\textbullet] $k$ is an empirical coefficient,
    \item[\textbullet] $\omega_{E\times B}$ is the flow shearing rate,
    \item[\textbullet] $\gamma_{ITG}$ is the ion temperature gradient growth rate,
    \item[\textbullet] $s_{thres}$ is the threshold shear,
    \item[\textbullet] $s$ is the magnetic shear.
\end{itemize}
This expression for $f_s$ was derived based on experimental results by \cite{pianroj2012simulations, sugihara2001simulation} to ensure that a transport barrier -the pedestal- is properly modelled near the plasma edge. \\ \\
Decoupling the magnetic effects from the thermal ones for our PDE-based control problem and defining $A = 2/(3a^2)$, we formualte the diffusion coefficient as
\begin{equation} \label{eq:diffusion}
    \chi(x,t) = A \bigl( B(x) + C(x) \sqrt{T_e(x,t)} \bigr) \bigl| \nabla T_e(x,t) \bigr|
\end{equation}
where
\begin{equation}
    \left\{
        \begin{array}{ll}
            B(x) &= \frac{8 \times10^{-5} R L_{T_e}}{B_{\phi0}} q^2(x) f_s(x) \\
            C(x) &= \frac{5 \times10^{-6}}{B_{\phi0}^2} f_s(x)
        \end{array}
    \right.
\end{equation}
and the constant $L_{T_e}$ measures the time-averaged temperature gradient from $\bigl(T_{e}(0.8)-T_{e}(1)\bigr)/T_{e}(1)$. \\ \\
In our model, the safety factor $q(x)$ is assumed to vary slowly compared to the electron temperature, justifying its time-decoupling from the thermal diffusivity formula.
\subsection{Optimal Control Problem} \label{subsec:OCP}
We formulate our problem as an optimal control problem in the Pontryagin framework with a continuum of horizons, the time-discretization of which simplifies to the more classical MPC-like receding horizon control. The goal is to reach a desired electronic temperature profile $\bar{T}_e(x)$ for the Tokamak plasma by controlling the net input heating power $u(x,t)$, while ensuring stability of the system.
The control law should not only minimize deviations from a dynamically evolving reference trajectory, but also be robust to errors due to inaccuracies in the model described by the dynamical equation
\begin{equation} \label{eq:state}
    \begin{split}
        \frac{\partial T_e}{\partial t} & = \frac{A}{x} \frac{\partial}{\partial x} \biggl( x \bigl( B(x) + C(x) \sqrt{T_e} \bigr) \Bigl( \frac{\partial T_e}{\partial x} \Bigr)^2 \biggr) + u \\
        & \stackrel{\eqref{eq:diffusion}}{=} \frac{1}{x} \frac{\partial}{\partial x} \Bigl( x \chi \frac{\partial T_e}{\partial x} \Bigr) + u
    \end{split}
\end{equation}
The optimal control problem then formulates as the PDE-constrained minimization problem of the distance to a predefined, dynamical trajectory leading to the final-time target at $t_f < \infty$, plus a regularity constraint on the $L^2_x = L^2(x dx)$-norm of the control variable
\begin{equation} \label{eq:ocp}
    \min_{\substack{T_e \in \mathcal{A} \\ u \in \mathcal{U}}} \biggl[ \int_0^{t_f} \Bigl( \mathcal{J}_1(t,u,T_e[u]) + \frac{\alpha(t)}{2} \| u \|_{L^2_x}^2 \Bigr) \, dt \biggr]
\end{equation}
where $\mathcal{A}$ is the set of admissible solutions to equation \eqref{eq:state}, and $\mathcal{U}$ is the admissible control space, assumed only to be $H^1$. The set of intermediate cost functionals $(\mathcal{J}_1(t,\cdot,\cdot))_{t \in [0,t_f]}$ and penalizing terms $(\alpha(t))_{t \in [0,t_f]}$ are detailed in the following sections.
\section{A Pontryagin-Based Approach} \label{sec:formalism}
\subsection{Receding Horizon Pontryagin's Principle} \label{subsec:PMP}
Pontryagin's classical setting with one end-point horizon has the inconvenience of providing an open-loop optimal controller independently from the real evolution of the system. It is therefore prone to errors due to numerical instabilities and deviations of the model from reality. Extending it through a smooth continuum of horizons could prevent these problems by providing a feedback optimal controller in a closed-loop with the system. Let the intermediate targets $( \hat{T}_e(\cdot,t))_{t \in [0,t_f]}$ define our receded horizons as an exponential interpolation between the initial state $T_{e,0}(\cdot)$ and the final-time target $\bar{T}_e(\cdot)$
\begin{equation}
    \hat{T}_e(x,t) = T_{e,0}(x) + (1 - e^{-\mu \nicefrac{t}{t_f}}) \bigl( \bar{T}_e(x) - T_{e,0}(x) \bigr)
\end{equation}
The corresponding intermediate cost functional measuring the distance of the controlled state to that intermediate target is
\begin{equation} \label{eq:costs}
    \mathcal{J}_1(t,u,T_e[u]) = \frac{1}{2} \int_0^1 \bigr( T_e(x,t) - \hat{T}_e(x,t) \bigl)^2 \, x dx
\end{equation}
Suppose that we controlled the system up to a time $t_i \in [0,t_f[$ along the reference trajectory $\hat{T}_e$ and set $t_{i+1} = t_i+\delta t \in ]t_i,t_f]$ a close enough time horizon. By this time discretization, the optimal control problem \eqref{eq:ocp} splits into subproblems as follows
\begin{equation}
    \sum_{\substack{0 \leq \cdot \cdot \leq t_i \leq \\ t_{i+1} \leq \cdot \cdot \leq t_f}} \delta t \min_{\substack{T_e \in \mathcal{A} \\ u \in \mathcal{U}}} \Bigl[ \mathcal{J}_1(t_{i+1},u,T_e[u]) + \frac{\alpha(t_{i+1})}{2} \| u \|_{L^2_x}^2 \Bigr]
\end{equation}
Without loss of generality, let us take any subinterval $[t,\tau]$ in the partition of $[0,t_f]$. Solving this subproblem of optimal control within Pontryagin's framework requires introducing an augmented Lagrangian $\mathcal{L}$ via the Lagrange multiplier $p$
\begin{align} \label{eq:Lagrangian}
    \mathcal{L}(\tau,u,T_e,p) & = \mathcal{J}_1(\tau,u,T_e) + \frac{\alpha(\tau)}{2} \int_t^\tau \int_0^1 u(x,t)^2 x \, dx + \biggl \langle p, \underbrace{\frac{1}{x} \frac{\partial}{\partial x} \Bigl( x \chi \frac{\partial T_e}{\partial x} \Bigr) + u - \frac{\partial T_e}{\partial t}}_{=0} \biggr \rangle_{L^2_x}
\end{align}
Separately calculating the PDE-constraint enforcing term via integrations by parts and using Green's formula gives
\begin{align}
    \biggl \langle p, \frac{1}{x} \frac{\partial}{\partial x} \Bigl( x \chi \frac{\partial T_e}{\partial x} \Bigr) + u - \frac{\partial T_e}{\partial t} \biggr \rangle_{L^2_x} & = \int_t^\tau \int_0^1 p u \, x dx dt + \int_t^\tau \cancel{x \chi p \frac{\partial T_e}{\partial x} \bigg|_0^1} \, dt \nonumber \\
    & \qquad - \int_t^\tau \int_0^1 \chi \frac{\partial p}{\partial x} \frac{\partial T_e}{\partial x} \, x dx dt + \int_t^\tau \int_0^1 \frac{\partial p}{\partial t} T_e \, x dx dt \\
    & \qquad - \int_0^1 p(x,\tau) T_e(x,\tau) \, x dx + \int_0^1 p(x,t) T_{e}(x,t) \, x dx \nonumber
\end{align}
The boundary term vanishes thanks to the homogeneous Neumann boundary condition $\frac{\partial T_e}{\partial x}(0,\cdot) = 0$ and by enforcing Dirichlet boundary condition $p(1,\cdot)=0$ on the costate. \\ \\
By the homogeneous Dirichlet boundary condition $T_e(1,\cdot) = 0$, the right hand-side's third term simplifies into
\begin{align}
    & - \int_t^\tau \int_0^1 \chi \frac{\partial p}{\partial x} \frac{\partial T_e}{\partial x} \, x dx dt = - \int_t^\tau \cancel{x \chi \frac{\partial p}{\partial x} T_e \bigg|_0^1} \, dt + \int_t^\tau \int_0^1 \frac{\partial}{\partial x} \Bigl( x \chi \frac{\partial p}{\partial x} \Bigr) T_e \, dx dt
\end{align}
Bringing together all the members of equation \eqref{eq:Lagrangian}, we get the final expression of the augmented Lagrangian
\begin{align}
    \mathcal{L}(\tau,u,T_e,p) & = \int_t^\tau \int_0^1 \Bigl[ p u + \frac{\alpha}{2} u^2 \Bigr] \, x dx dt + \int_0^1 \Bigl[ \frac{1}{2} \bigr( T_e(x,\tau) - \hat{T}_e(x,\tau) \bigl)^2 - p(x,\tau) T_e(x,\tau) \Bigr] \, x dx \nonumber \\
    & \qquad + \int_t^\tau \int_0^1 \biggl[ \frac{\partial p}{\partial t} + \frac{1}{x} \frac{\partial}{\partial x} \Bigl( x \chi \frac{\partial p}{\partial x} \Bigr) \biggr] T_e \, x dx dt + \int_0^1 p(x,t) T_{e}(x,t) \, x dx
\end{align}
\begin{remark}
    The last term is irrelevant to our optimization, since it no longer is a horizon but a given "initial condition" for the optimal control problem on the interval $[t,\tau]$.
\end{remark}
For simplicity, assume that the diffusivity $\chi(T_e) = \chi$ is state-independent on $[0,1] \times [t, \tau]$ (since $B(x), C(x) \ll 1$).
Set the Lagrangian stationary with respect to the state
\begin{equation}
    0 =
    \begin{pmatrix}
        \nabla_{T_e(\cdot,\cdot)} \mathcal{L} \\
        \nabla_{T_{e}(\cdot,\tau)} \mathcal{L}
    \end{pmatrix}
\end{equation}
the latter condition rewrites as the adjoint equation
\begin{equation} \label{eq:costate}
    \left\{
        \begin{array}{ll}
            \frac{\partial p}{\partial t} = - \frac{1}{x} \frac{\partial}{\partial x} \Bigl( x \chi \frac{\partial p}{\partial x} \Bigr), \quad \text{a.e. on } [0,1] \times [t, \tau] \\
            p(x,\tau) = T_e(x,\tau) - \hat{T}_e(x,\tau), \quad \text{a.e. in } [0,1]
        \end{array}
    \right.
\end{equation}
with Dirichlet boundary condition $p(1,\cdot) = 0$. Making the Lagrangian stationary with respect to the control yields Pontryagin's optimality condition
\begin{equation} \label{eq:optimality}
    0 = \nabla_u \mathcal{L}(\tau) = p + \alpha(\tau) u
\end{equation}
Injecting this expression back into the state equation \eqref{eq:state} yields the coupled state-costate PDE system on $[0,1] \times [t, \tau]$
\begin{equation} \label{eq:coupled}
    \left\{
        \begin{array}{ll}
            \frac{\partial T_e}{\partial t} = \frac{1}{x} \frac{\partial}{\partial x} \Bigl( x \chi \frac{\partial T_e}{\partial x} \Bigr) - \alpha^{-1} p \\
            \frac{\partial p}{\partial t} = - \frac{1}{x} \frac{\partial}{\partial x} \Bigl( x \chi \frac{\partial p}{\partial x} \Bigr)
        \end{array}
    \right.
\end{equation}
and the following initial-final conditions
\begin{equation}
    \left\{
        \begin{array}{ll}
            T_e(\cdot,t) = T_{e,t}(\cdot) \\
            p(\cdot,\tau) = T_e(\cdot,\tau) - \hat{T}_e(\cdot,\tau)
        \end{array}
    \right.
\end{equation}
Except for the transversality conditions, notice that PMP does not depend on the time horizon choice, so we differentiate the former by making $\delta t$ arbitrarily small
\begin{equation}
    \begin{split}
        \frac{\partial p}{\partial t} & = \lim_{\delta t \to 0} \frac{p(\cdot,t+\delta t) - p(\cdot,t)}{\delta t} \\
        & = \lim_{\tau \downarrow t} \frac{ \bigl( T_e(\cdot,\tau) - \hat{T}_e(\cdot,\tau) \bigr) - \bigl( T_e(\cdot,t) - \hat{T}_e(\cdot,t) \bigr)}{\tau - t} \\
        & = \frac{\partial T_e}{\partial t} - \frac{\partial \hat{T}_e}{\partial t} 
    \end{split}
\end{equation}
Substituting back into the state-costate equations \label{eq:cmd}
\begin{equation} 
    - \frac{1}{x} \frac{\partial}{\partial x} \Bigl( x \chi \frac{\partial p}{\partial x} \Bigr) = \frac{1}{x} \frac{\partial}{\partial x} \Bigl( x \chi \frac{\partial T_e}{\partial x} \Bigr) - \frac{1}{\alpha} p - \frac{\partial \hat{T}_e}{\partial t}
\end{equation}
Hence the inverse problem on the adjoint variable $p$
\begin{equation}
    \biggl( \frac{1}{\alpha} Id - \frac{1}{x} \frac{\partial}{\partial x} \Bigl( x \chi \frac{\partial}{\partial x} \Bigr) \biggr) p = \frac{1}{x} \frac{\partial}{\partial x} \Bigl( x \chi \frac{\partial T_e}{\partial x} \Bigr) - \frac{\partial \hat{T}_e}{\partial t}
\end{equation}
Combining both transversality and stationarity conditions into this quasi-steady solution for the adjoint state, we bypass the difficulty of solving it backward in time as in the original one-end horizon formulation. \\ \\
We derived a closed-loop optimal controller evolving forward in time alongside the controlled state, guiding the system along a reference trajectory toward the final target and enabling real-time error corrections by feedback.
\subsection{Regularity Analysis} \label{subsec:reg}
In this section, we analyse the regularity of the adjoint state $p$ and thus that of the optimal controller $u$ via Pontryagin's optimality condition \eqref{eq:optimality} in the sense of a certain $W^{m,2}_x = H^m_x$ Sobolev space ($m=1$ in the sequel). \\ \\
Multiplying \eqref{eq:cmd} by $p$ and $L^2_x$-integrating over the domain
\begin{equation}
    \begin{split}
        & \int_0^1 \alpha^{-1} p^2 \, x dx - \int_0^1 \frac{1}{x} \frac{\partial}{\partial x} \Bigl( x \chi \frac{\partial p}{\partial x} \Bigr) p \, x dx = \int_0^1 \frac{1}{x} \frac{\partial}{\partial x} \Bigl( x \chi \frac{\partial T_e}{\partial x} \Bigr) p \, x dx - \int_0^1 \frac{\partial \hat{T}_e}{\partial t} p \, x dx
    \end{split}
\end{equation}
Using Green's identity and the boundary conditions
\begin{equation}
    \begin{split}
        & \int_0^1 \alpha^{-1} p^2 \, x dx + \int_0^1 \chi \frac{\partial p}{\partial x}^2 x dx - \cancel{x \chi p \frac{\partial p}{\partial x} \bigg|_0^1} = \cancel{x \chi p \frac{\partial T_e}{\partial x} \bigg|_0^1} - \int_0^1 \chi \frac{\partial T_e}{\partial x} \frac{\partial p}{\partial x} \, x dx - \int_0^1 \frac{\partial \hat{T}_e}{\partial t} p \, x dx
    \end{split}
\end{equation}
Cauchy-Schwarz inequality applies to the right hand-side
\begin{flalign} \label{eq:ineq}
    \int_0^1 \alpha^{-1} p^2 \, x dx + \int_0^1 \chi \frac{\partial p}{\partial x}^2 x dx & \leq \sqrt{\int_0^1 \chi \frac{\partial T_e}{\partial x}^2 x dx} \sqrt{\int_0^1 \chi \frac{\partial p}{\partial x}^2 x dx} \nonumber \\
    & \qquad \qquad + \sqrt{\int_0^1 \frac{\partial \hat{T}_e}{\partial t}^2 x dx} \sqrt{\int_0^1 p^2 x dx}
\end{flalign}
Young's inequality applies to the right hand-side's terms
\begin{equation}
    \begin{split}
       \sqrt{\int_0^1 \chi \frac{\partial T_e}{\partial x}^2 x dx} & \sqrt{\int_0^1 \chi \frac{\partial p}{\partial x}^2 x dx} \leq \frac{1}{2} \int_0^1 \chi \frac{\partial T_e}{\partial x}^2 x dx + \frac{1}{2} \int_0^1 \chi \frac{\partial p}{\partial x}^2 x dx 
    \end{split}
\end{equation}
along its positive $\varepsilon$-tradeoff Peter-Paul inequality version
\begin{equation}
    \begin{split}
        \sqrt{\int_0^1 \frac{\partial \hat{T}_e}{\partial t}^2 x dx} & \sqrt{\int_0^1 p^2 x dx} \leq \frac{1}{4 \varepsilon} \int_0^1 \frac{\partial \hat{T}_e}{\partial t}^2 x dx + \varepsilon \int_0^1 p^2 x dx
    \end{split}
\end{equation}
Plugging these new upper bounds back into \eqref{eq:ineq} yields
\begin{equation}
    2 (\alpha^{-1} - \varepsilon) \| p \|_{L^2_x}^2 + \| p \|_{\dot{H}^1_{x \chi}}^2 \leq \| T_e \|_{\dot{H}^1_{x \chi}}^2 + \frac{1}{2 \varepsilon} \| \partial_t \hat{T}_e \|_{L^2_x}^2 
\end{equation}
Maintaining coercivity of the left hand-side's norms requires choosing $0 < \varepsilon \leq \alpha^{-1}$ in the sequel.
\begin{theorem}[Weighted Poincaré Inequality]
    Let $w \in L^\infty(\Omega)$ be a nonnegative function such that $w(x) > 0$ a.e.. Then there exists a constant $C > 0$ such that for all $v \in H^1_w(\Omega)$
    \begin{equation}
        \int_\Omega w(x) \, v(x)^2 \, dx \leq C \int_\Omega w(x) |\nabla v(x)|^2 \, dx
    \end{equation}
    The optimal Poincaré constant is $C = \lambda_{1,w}(\Omega)^{-1}$, where $\lambda_{1,w}$ is the first Dirichlet-Laplacian eigenvalue of $\Omega$ as the infimum of the weighted Rayleigh quotient $\mathcal{R}_w$ over $H^1_w(\Omega)$
    \begin{equation}
        \lambda_{1,w}(\Omega) = \inf_{v \in H^1_w(\Omega)} \bigl[ \mathcal{R}_w(v) := \| v \|_{\dot{H}_w^1(\Omega)}^2 / \| v \|_{L_w^2(\Omega)}^2 \bigr]
    \end{equation}
    and $\| \cdot \|_{\dot{H}_w^1(\Omega)}$ is the $w$-weighted Sobolev seminorm defined for all $v \in H^1_w(\Omega)$ by $\| v \|_{\dot{H}_w^1(\Omega)}^2 = \int_\Omega |\nabla v|^2 \, w dx$.
\end{theorem}
Poincaré inequality holds in this case since the nonnegative weight $w = x \chi \in L^\infty(0,1)$ is almost everywhere positive
\begin{equation}
    \lambda_{1,x \chi} \| p \|_{L^2_{x \chi}}^2 \leq \| p \|_{\dot{H}^1_{x \chi}}^2
\end{equation}
Conveniently choosing $\varepsilon = \frac{1}{2 \alpha}$ for illustrative purposes and combining the inequalities above gives one final estimate
\begin{equation}
    \alpha^{-1} \| p \|_{L^2_x}^2 + \lambda_{1,x \chi} \| p \|_{L^2_{x \chi}}^2 \leq \| T_e \|_{\dot{H}^1_{x \chi}}^2 + \alpha \| \partial_t \hat{T}_e \|_{L^2_x}^2 
\end{equation}
which translates into an upper bound for the energy injected via the optimal controller using equation \eqref{eq:optimality}
\begin{equation} \label{eq:bound}
    \| u \|_{L^2_x}^2 + \lambda_{1,x \chi} \alpha \| u \|_{L^2_{x \chi}}^2 \leq \| \partial_t \hat{T}_e \|_{L^2_x}^2 + \alpha^{-1} \| T_e \|_{\dot{H}^1_{x \chi}}^2 
\end{equation}
\begin{remark}
    The first term on the right hand-side is finite by construction, since the reference $\hat{T}_e \in H^1((0,t_f),L^2_x)$ was chosen as a regular interpolation so that $\partial_t \hat{T}_e \in L^2_x$, and the second term is proportional to the Dirichlet energy of a controlled diffusion equation in finite time.
\end{remark}
From the inequalities above and the finiteness of the upper-bound in inequality \eqref{eq:bound}, it follows that the control variable $u$ and the adjoint state $p$ are both regular in the sense of $H^1_{x \chi}$ Sobolev space and bounded in the $L^2_x$-energy space. \\ \\
As a measure of the total energy injected into the system, we naturally chose the $L^2_x$-energy norm of the control variable $u$ as bounded in inequality \eqref{eq:bound} (recall $\lambda_1 > 0)$
\begin{equation} \label{eq:estimate}
    \| u \|_{L^2_x} \leq \sqrt{\| \partial_t \hat{T}_e \|_{L^2_x}^2 + \alpha^{-1} \| T_e \|_{\dot{H}^1_{x \chi}}^2}
\end{equation}
For future refined analysis of the boundedness/regularity, one would need to study the time-evolution of the upper-bound in inequality \eqref{eq:bound}. Here below is one way to do such
\begin{lemma}[Gradient Flow]
    The free diffusion equation is the descending gradient flow of the Dirichlet energy.
\end{lemma}
\begin{align}
    \frac{\partial}{\partial t} \| T_e \|_{\dot{H}^1_\chi(\Omega)}^2 & = \Bigl \langle \nabla_{T_e} \| \cdot  \|_{\dot{H}^1_\chi(\Omega)}^2 , \frac{\partial}{\partial t} T_e \Bigr \rangle_{L^2(\Omega)} \\
    & = \Bigl \langle - \nabla \cdot ( \chi \nabla T_e) , \nabla \cdot ( \chi \nabla T_e) \Bigr \rangle_{L^2(\Omega)} \leq 0 \nonumber
\end{align}
However, we are dealing with a diffusion equation that is controlled rather than freely evolving, hence the need to add back our control variable $u$ in the last bracket term and perform adequate functional-analytic manipulations on its norm nested in both hand-sides of inequality $\eqref{eq:bound}$. We settle with reporting results in figure (\ref{fig:bound}) in the sequel.
\subsection{Convergence Analysis} \label{subsec:conv}
We want to prove convergence of the controlled state toward the final target, which is measured by the objective functional
\begin{equation}
    \mathcal{J}(t) =  \frac{1}{2} \int_0^1 \bigr( T_e(x,t) - \bar{T}_e(x) \bigl)^2 \, x dx
\end{equation}
Inserting the intermediate target $\hat{T}_e(\cdot,t)$ between both terms and using the definition of $\hat{T}_e(\cdot,t_f)$
\begin{equation} \label{eq:J}
    \begin{split}
        \mathcal{J}(t) & = \frac{1}{2} \int_0^1 \Bigr( \bigl( T_e(x,t) - \hat{T}_e(x,t) \bigr) + \bigl( \hat{T}_e(x,t) - \hat{T}_e(x,t_f) \bigr) \Bigl)^2 \, x dx
    \end{split}
\end{equation}
By virtue of the classical arithmetic–geometric inequality
\begin{equation} \label{eq:J12}
    \mathcal{J}(t) \leq 2 \bigl( \mathcal{J}_1(t) + \mathcal{J}_2(t) \bigr)
\end{equation}
where the two functionals $\mathcal{J}_1$ and $\mathcal{J}_2$ are defined as
\begin{equation}
    \left\{
        \begin{array}{ll}
            \mathcal{J}_1(t) & = \frac{1}{2} \int_0^1 \bigl( T_e(x,t) - \hat{T}_e(x,t) \bigr)^2 \, x dx \\
            \mathcal{J}_2(t) & = \frac{1}{2} \int_0^1 \bigl( \hat{T}_e(x,t) - \hat{T}_e(x,t_f) \bigr)^2 \, x dx
        \end{array}
    \right.
\end{equation}
Trivially, we see that the second term $\mathcal{J}_2$ is exponentially decaying. Indeed, plugging in the expression of the reference trajectory $\hat{T}_e(x,t)$ as an exponential interpolation gives
\begin{equation}
    \begin{split}
        \mathcal{J}_2(t) & = \frac{1}{2} \int_0^1 \bigl( \hat{T}_e(x,t) - \hat{T}_e(x,t_f) \bigr)^2 \, x dx \\
        & = \frac{1}{2} \int_0^1 e^{-2 \mu \nicefrac{t}{t_f}} \bigl( \bar{T}_e(x,0) - \hat{T}_e(x,t_f) \bigr)^2 \, x dx \\
        & = e^{-2 \mu \nicefrac{t}{t_f}} \mathcal{J}_2(0), \quad t \in [0,t_f]
    \end{split}
\end{equation}
As for the first term $\mathcal{J}_1$, we will show that it remains small according to our optimization-based control method, aiming at minimizing the deviation of the controlled state from a given reference trajectory, while calibrating in real-time the parameter $\alpha(t)$ based on the evolution of $\mathcal{J}_1(t)$ itself.

More specifically, we want to establish an empirical law governing the penalty term $\alpha(t)$ of the form
\begin{equation}
    \dot{\alpha}(t) = \beta(t) \times \, h(\mathcal{J}_1(t),\dot{\mathcal{J}}_1(t))
\end{equation}
where $\beta \in \{-1,1\}$ is the adjustment sign to compensate overshooting or undershooting a midway-target, defined as
\begin{equation} \label{eq:sign}
    \beta(t) = \text{sign} \biggl( \int_0^1 \bigr( T_e(x,t) - \hat{T}_e(x,t) \bigl) \, x dx \biggr)
\end{equation}
and $h \geq 0$ measures the strength of the deviation of $\mathcal{J}_1$ from $0$ as detailed below. Intuitively, if $T_e(\cdot,t)$ is above the midway-target $\hat{T}_e(\cdot,t)$ on average , the sign is positive and we amplify the regularity constraint on the control, otherwise we relax it. \\

To study the effect of $\alpha$, we intend to calibrate it once using the standard one-horizon Pontryagin setting on the whole time interval $[0,t_f]$ instead of $[t,\tau]$ as developed in Section \ref{subsec:PMP}. We retrieve the same adjoint equation \eqref{eq:costate} but on $[0,1] \times [0,t_f]$ alongside Pontryagin's optimality conditions \eqref{eq:optimality} leading to the coupled state-costate PDEs \eqref{eq:coupled} with initial conditions at time $t=0$ and transversality conditions at final time $\tau=t_f$.

Solving the above coupled PDEs using algorithm (\ref{algo}) below for a constant penalty term $\alpha$ gives us the optimal energy $\mathcal{J}^*(\alpha) = \mathcal{J}(u^*,T_e^*,p^*)$ at the optimal triplet $(u^*,T_e^*,p^*)$ solution to PMP with $\alpha$-regularized Lagrangian $\mathcal{L}$. Doing the same for different values of $\alpha$ allows us to heuristically analyse the variations of $\mathcal{J}^*$ with respect to $\alpha$ around the optimal penalty term $\alpha^* = \arg \min_{\alpha} \mathcal{J}^*(\alpha)$.
\begin{remark}
    Choosing a small penalty term $\alpha$ leads to a solution overshooting the target and vice-versa; a large $\alpha$ leads to a solution undershooting the target.
\end{remark}
The following quadratic relation best approximates the results obtained from numerical experimentations
\begin{equation}
    \mathcal{J}^*(\alpha) = \kappa (\alpha - \alpha^*)^2, \quad \kappa > 0
\end{equation}
This empirical relation is valid for a classical one-horizon Pontryagin framework at a fixed time. Extending it to our continuum of horizons setting requires further assumptions on the coefficient $\kappa$ and the evolving optimal regularizing parameter, that is their time-invariability $\kappa(t)=\kappa$ and $\alpha^*(t)=\alpha^*$
\begin{equation} \label{eq:param}
    \mathcal{J}_1(\alpha(t)) = \kappa (\alpha(t) - \alpha^*)^2
\end{equation}
Under these hypothesises, we derive the empirical formula
\begin{equation}
    \begin{split}
        \dot{\mathcal{J}}_1(t) = \frac{d}{dt} \mathcal{J}_1(\alpha(t)) & = \frac{\partial \mathcal{J}_1}{\partial \alpha} \cdot \dot{\alpha}(t) \\
        & \stackrel{\eqref{eq:param}}{=} 2 \kappa (\alpha(t) - \alpha^*) \cdot \dot{\alpha}(t) \\
        & \stackrel{\eqref{eq:sign}}{=} 2 \beta \sqrt{\kappa} \sqrt{\mathcal{J}_1(t)} \cdot \dot{\alpha}(t) \\ 
        \Rightarrow |\dot{\alpha}(t)| = h(\mathcal{J}_1,\dot{\mathcal{J}}_1) & \propto \frac{\dot{\mathcal{J}}_1(t)}{\sqrt{\mathcal{J}_1(t)}} 
    \end{split}
\end{equation}
By including this adaptive feedback mechanism for the regularizing parameter in our algorithm, we complete the synthesis of the optimal control law.
\section{Control Algorithm Implementation and Simulation Results} \label{sec:results}
\subsection{Control Algorithm} \label{subsec:algo}
The proposed adjoint-based optimal control algorithm, accessible in Subsection \ref{subsec:code}, is specifically tuned to conform with the H-mode configuration of Tore Supra Tokamak for tracking a desired electronic temperature profile $\bar{T}_e$. Simulation results are presented in Subsection \ref{subsec:simul} below. Hereafter is a pseudocode for solving open-loop PMP, which was independently used to fine-tune $\alpha^* = 10$ from experimentation as studied at the end of Subsection \ref{subsec:conv}.
\begin{algorithm}[H]
\caption{\centering Forward-Backward Sweep} \label{algo}
\textbf{Input} $\alpha, \varepsilon=10^{-6}, n_{max}=10^3$ \\
\textbf{Output} $u,y,p$
    \begin{algorithmic}
        \State
            $n \gets 0$
        \State 
            $u_0 \gets 0$
        \While{$\|u_{n}-u_{n-1}\|_{L^2([0,t_f],L^2(\Omega))} > \varepsilon$ and $n < n_{max}$}
            \State
                $y_n \gets \text{solution to the } u_n \text{-controlled state equation}$  
            \State 
                $p_n \gets \text{solution to the } y_n \text{-conditioned adjoint equation}$
            \State 
                $u_n \gets - \alpha^{-1} \nabla_u \mathcal{L}(u_n,y_n,p_n)$
            \State
                $n \gets n+1$
        \EndWhile
        \State 
            $u \gets u_n, \; y \gets y_n, \; p \gets p_n$
    \end{algorithmic}
\end{algorithm}
\subsection{Simulation Results} \label{subsec:simul}
\begin{figure}[H]
    \centering
    \includegraphics[width=0.9\textwidth]{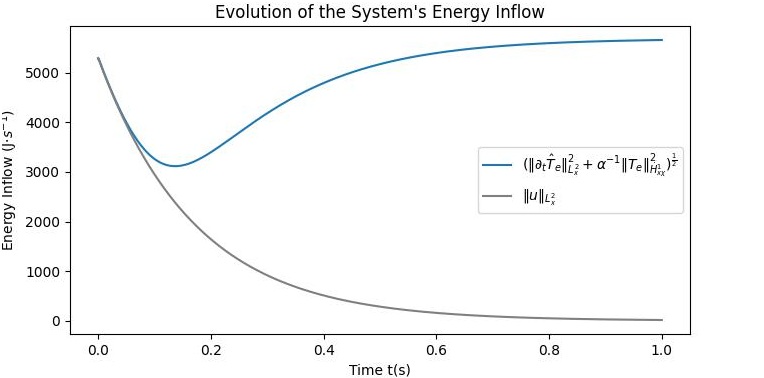}
    \caption{\centering Boundedness of the Optimal Power Input}
    \label{fig:bound}
\end{figure}

\begin{figure}[H]
    \centering
    \includegraphics[width=0.9\textwidth]{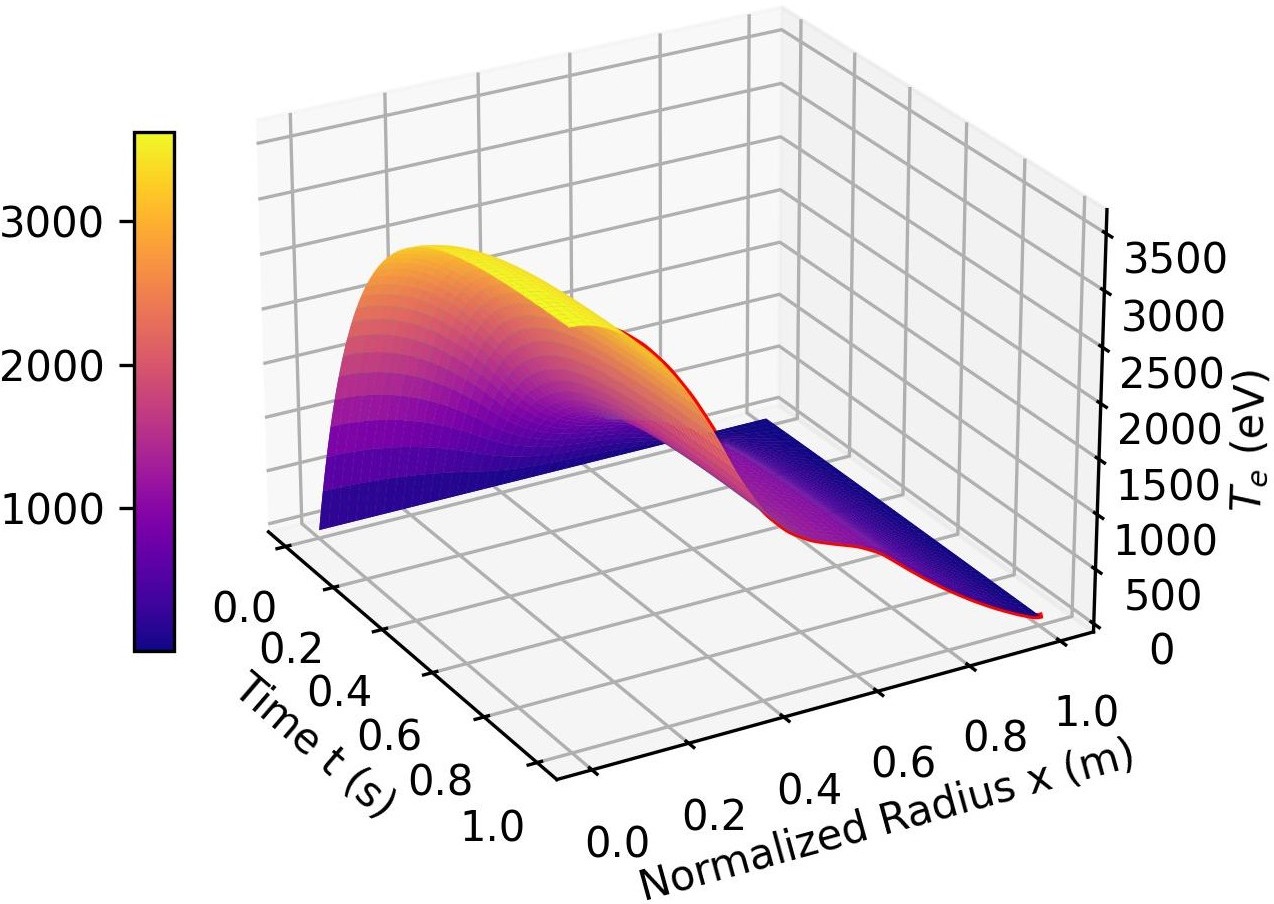}
    \caption{\centering Controlled Plasma Temperature Evolution}
    \label{fig:state}
\end{figure}

\begin{figure}[H]
    \centering
    \includegraphics[width=0.9\textwidth]{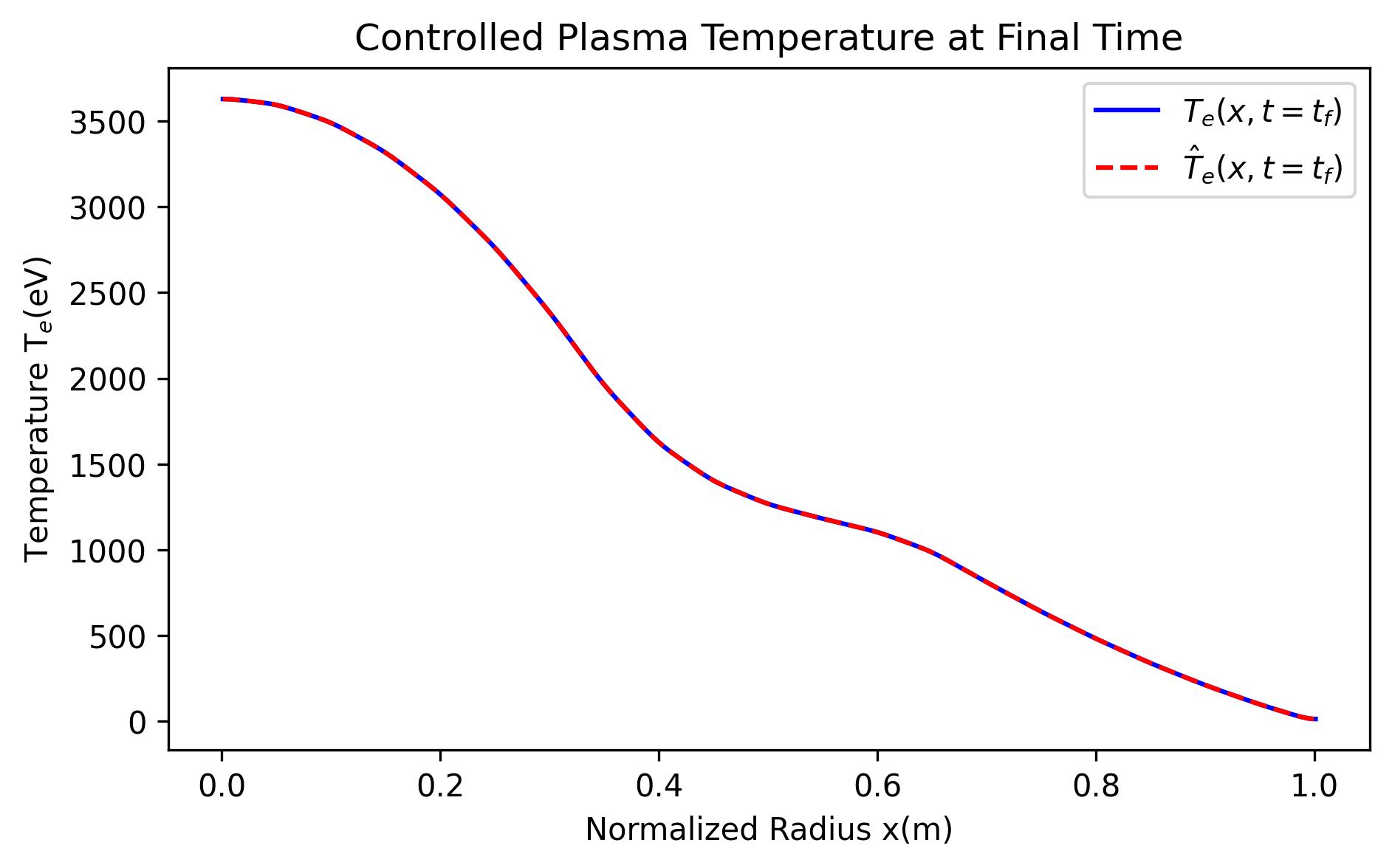}
    \caption{\centering Controlled Plasma Temperature at Final Time}
    \label{fig:target}
\end{figure}

\begin{figure}[H]
    \centering
    \includegraphics[width=0.9\textwidth]{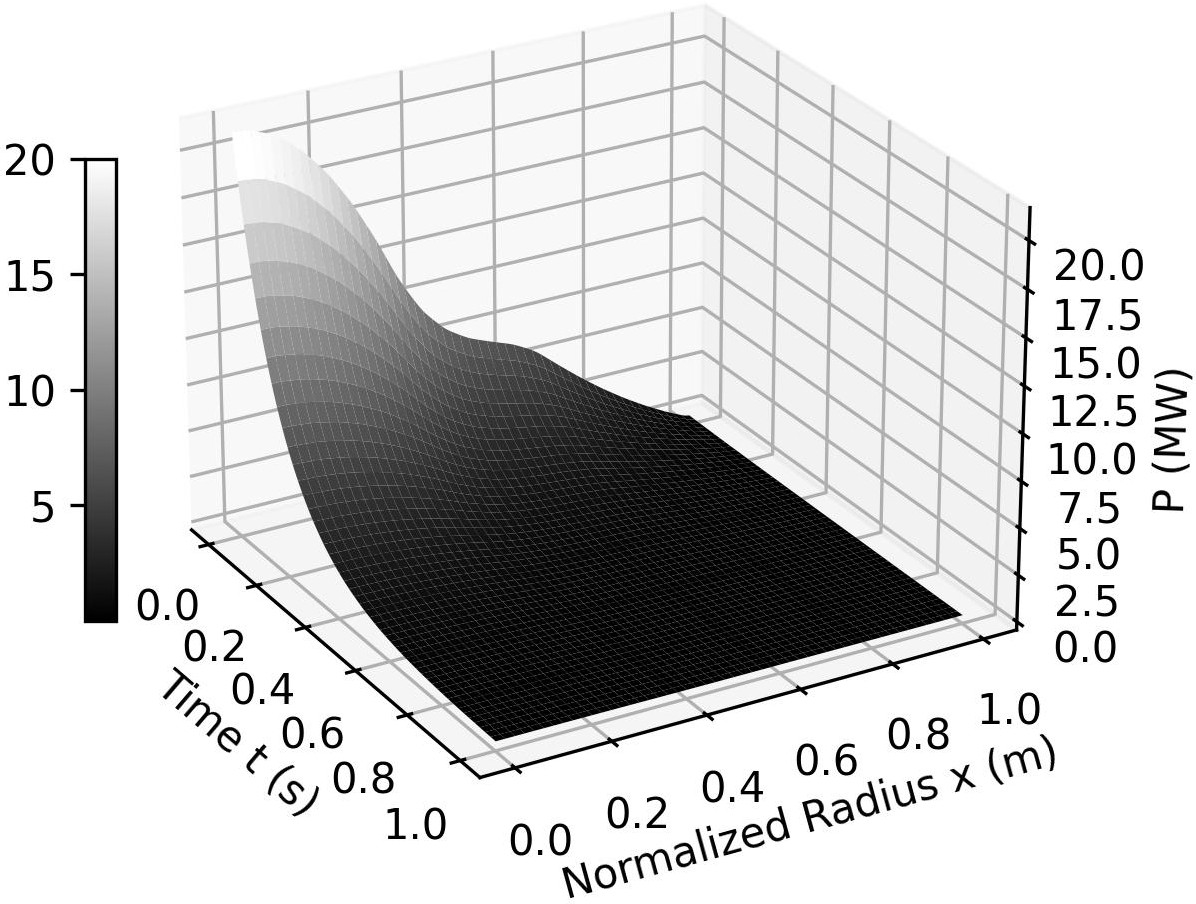}
    \caption{\centering Optimal Power Input Control Law}
    \label{fig:control}
\end{figure}

\begin{figure}[H]
    \centering
    \includegraphics[width=0.9\textwidth]{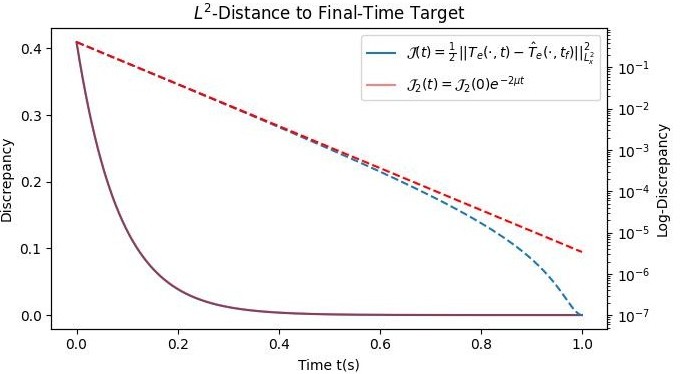}
    \caption{\centering Convergence of the Objective Functional}
    \label{fig:convergence}
\end{figure}

All numerical simulations are performed by taking into account the model presented in Subsection \ref{subsec:model}, where the parameters of the heat diffusivity are calibrated on experimental measurements extracted from the Tore Supra shot 36 056 (12/08/2005) obtained with lower hybrid (2.5 MW) and ion cyclotron (5 MW) radio frequency antennas. The free parameters used in Subsection \ref{subsec:PMP} are as follows: the final time is set to $t_f=1$, and the exponentiation parameter in the expression of the reference trajectory is fixed at $\mu=5.85$ for an abrupt ramp-up phase at the start of the plasma heating. \\ \\
Figure \ref{fig:bound} reports the time-evolution of the total mean energy inflow needed for the system's control toward an optimal functioning regime as analysed in Subsection \ref{subsec:reg}. The grey curve represents the $L^2$-energy norm of the control, while the blue curve corresponds to the theoretical upper bound derived from functional regularity estimates.
From an operational standpoint, maintaining a tight upper bound on the control energy during the ramp-up phase is particularly important for the system's transition toward an optimal confinement regime with minimal expenditure. This energetic efficiency is crucial for setting down power balances given the expected energy release from plasma fusion reactions in such high confinement modes. Despite this upper bound getting loose with time, based on the rough estimate \eqref{eq:estimate}, most of the heating energy is rather spent on the fusion plasma ignition during the first phase than on the plasma stabilization. \\ \\
Figure \ref{fig:state} depicts the space-time evolution of the controlled plasma temperature $T_e(x,t)$ along the prescribed reference trajectory $\hat{T}_c(x,t)$, which was defined as an exponential interpolation connecting the initial state to the experimentally identified stationary target corresponding to the Tore Supra optimal H-mode equilibrium. 
The controlled dynamics remain nearly indistinguishable from the reference trajectory, demonstrating the high fidelity of the control strategy in tracking the desired temperature evolution. Figure \ref{fig:target} further confirms that the controlled temperature profile almost-perfectly matches the final-time target with a discrepancy of $1.01 \cdot 10^{-7}$. 
The space–time distribution of the optimal power input in Figure \ref{fig:control} clearly illustrates the physically expected shape of auxiliary heating in an H-mode Tokamak discharge: an initial, sharp energy deposition localized near the plasma core efficiently triggers the transition to high confinement, while the required power rapidly decreases as the temperature gradient self-sustains through improved confinement. This behavior reproduces the typical pattern observed in high-performance discharges, where feedback actuators operate predominantly during the transient phase.
Experimentally, the evolution of the optimal controller both reaffirms the energetic feasibility of the control process and validates the estimate derived in Subsection \ref{subsec:reg} regarding the boundedness of the control variable. \\ \\
As studied in Subsection \ref{subsec:conv}, the convergence analysis is further supported by the plots in Figure \ref{fig:convergence}: two almost-overlapping curves; that of the exponentially decaying reference trajectory in red and that of the main objective functional in blue nearly fitting the former up to a small deviation due to our optimization-based control technique, as rendered more clearly visible on their log-scaled dotted counterparts.
The exponential speed of convergence is essentially the same as that imposed on the controlled system through the reference trajectory, until orders of magnitudes below where tiny deviations play out due to control-induced minimization of the first bounding term in inequality \eqref{eq:J12}. \\ \\
Figure \ref{fig:convergence-bis} below showcases the evolution of $J_1$ as a measure of deviation from predefined intermediate targets, topping out at $3.369 \cdot 10^{-6}$ by the end of the control process
\begin{figure}[H]
    \centering
    \includegraphics[width=0.9\textwidth]{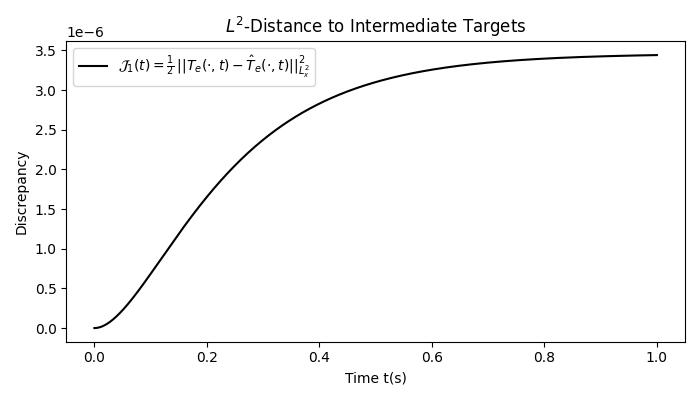}
    \caption{\centering Tracking the Reference Trajectory}
    \label{fig:convergence-bis}
\end{figure}
These deviations have been controlled to be as small as possibly allowed by our regularized optimization-based control within Pontryagin's framework, while blending the adaptive law of the regularizating $\alpha$-parameter into our control algorithm.
\subsection{Code Script} \label{subsec:code}
The implemented control algorithm is available at:  
\url{https://gitlab.com/slimj-group/SlimJ-project/-/blob/main/Optimal_Control_of_Tokamak_Plasma_Temperature_via_Pontryagin_Receding_Horizon_Principle.ipynb}
\section{Conclusion} \label{sec:concl}
In this work, we introduced a Pontryagin-framed optimal control formulation tailored to the regulation of temperature profiles in H-mode Tokamak plasmas. By extending the classical single-horizon Pontryagin framework into a continuum of receding intermediate targets, the proposed approach merges the analytical rigor PMP with the practical efficiency of RHC. This hybrid formulation effectively transforms an open-loop optimal controller into a closed-loop feedback mechanism, enabling real-time compensation of model uncertainties and enhanced robustness to external perturbations.

We demonstrated that the resulting control law ensures exponential convergence of the plasma temperature toward its prescribed optimal profile, while satisfying physical constraints of Tokamak plasmas. Numerical validations under the Tore Supra configuration confirmed both the accuracy and energetic efficiency of the proposed scheme compared to conventional methods.

Beyond its computational simplicity relative to Lyapunov-based approaches in infinite-dimensional settings, the present method establishes a foundational bridge between optimal control theory and plasma physics applications. Furthermore, it opens a pathway toward the control of coupled nonlinear PDEs governing magnetically confined fusion plasmas, where magnetohydrodynamic phenomena at all scales may strongly interact within the same system.

Future investigations will extend this framework to incorporate geometric shape state constraints and explore its real-time implementation on TCV’s RAPTOR simulator, with a view toward integration in reactor-scale devices such as ITER. In this perspective, the proposed method applies to different configurations of Tokamak plasmas, offering scalable control synthesis to achieve robust, energy-efficient plasma fusion powerplants for future generations.
\bibliography{ifacconf}
\end{document}